\documentclass[12pt]{amsart}
\usepackage{amssymb}
\usepackage{color}
\usepackage{amsmath,epic,curves,amscd}
\usepackage[english]{babel}
\usepackage{graphicx}
\usepackage{comment}
\usepackage{appendix}
\usepackage{mathdots}
\usepackage[all]{xy}
\theoremstyle{remark}

\begin{document}
\baselineskip 6.0 truemm
\parindent 1.5 true pc

\newcommand\lan{\langle}
\newcommand\ran{\rangle}
\newcommand\tr{{\text{\rm Tr}}\,}
\newcommand\ot{\otimes}
\newcommand\ol{\overline}
\newcommand\join{\vee}
\newcommand\meet{\wedge}
\renewcommand\ker{{\text{\rm Ker}}\,}
\newcommand\image{{\text{\rm Im}}\,}
\newcommand\id{{\text{\rm id}}}
\newcommand\tp{{\text{\rm tp}}}
\newcommand\pr{\prime}
\newcommand\e{\epsilon}
\newcommand\la{\lambda}
\newcommand\inte{{\text{\rm int}}\,}
\newcommand\ttt{{\text{\rm t}}}
\newcommand\spa{{\text{\rm span}}\,}
\newcommand\conv{{\text{\rm conv}}\,}
\newcommand\rank{\ {\text{\rm rank of}}\ }
\newcommand\re{{\text{\rm Re}}\,}
\newcommand\ppt{\mathbb T}
\newcommand\rk{{\text{\rm rank}}\,}
\newcommand\SN{{\text{\rm SN}}\,}
\newcommand\SR{{\text{\rm SR}}\,}
\newcommand\HA{{\mathcal H}_A}
\newcommand\HB{{\mathcal H}_B}
\newcommand\HC{{\mathcal H}_C}
\newcommand\CI{{\mathcal I}}
\newcommand{\bra}[1]{\langle{#1}|}
\newcommand{\ket}[1]{|{#1}\rangle}
\newcommand\cl{\mathcal}
\newcommand\idd{{\text{\rm id}}}
\newcommand\OMAX{{\text{\rm OMAX}}}
\newcommand\OMIN{{\text{\rm OMIN}}}
\newcommand\diag{{\text{\rm Diag}}\,}
\newcommand\calI{{\mathcal I}}
\newcommand\bfi{{\bf i}}
\newcommand\bfj{{\bf j}}
\newcommand\bfk{{\bf k}}
\newcommand\bfl{{\bf l}}
\newcommand\bfp{{\bf p}}
\newcommand\bfq{{\bf q}}
\newcommand\bfzero{{\bf 0}}
\newcommand\bfone{{\bf 1}}
\newcommand\im{{\mathcal R}}
\newcommand\ha{{\frac 12}}
\newcommand\xx{{\text{\sf X}}}

\title{Indecomposable exposed positive bi-linear
maps between two by two matrices}

\author{Seung-Hyeok Kye}
\address{Department of Mathematics and Institute of Mathematics, Seoul National University, Seoul 151-742, Korea}
\email{kye at snu.ac.kr}
\thanks{partially supported by NRF-2017R1A2B4006655.}



\begin{abstract}
Positive bi-linear maps between matrix algebras play important roles
to detect tri-partite entanglement by the duality between bi-linear
maps and tri-tensor products. We exhibit indecomposable positive
bi-linear maps between $2\times 2$ matrices which generate
extreme rays in the cone of all positive bi-linear maps. In fact, they
are exposed, and so detect entanglement of positive partial transpose whose volume is nonzero.
\end{abstract}

\maketitle

\section{Introduction}

Positive linear maps in matrix algebras  are
indispensable to distinguish entanglement from separability, by the
duality \cite{{eom-kye},{horo-1}} between bi-partite separable
states and positive linear maps. The convex cone consisting of
positive linear maps is very complicated even in the cases when the
domains and ranges are low dimensional matrices, and its whole
structures are far from to be understood. Typical examples of
positive maps from $m\times m$ matrices into $n\times n$ matrices
are of the forms
\begin{equation}\label{dec-exp}
X\mapsto V^*XV,\qquad X\mapsto W^*X^\ttt W,
\end{equation}
with the transpose  $X^\ttt$ of $X$ and $m\times n$ matrices $V$ and
$W$. A sum of the first (respectively second) type is called
completely positive (respectively completely copositive), and a sum
of them is said to be decomposable. It is known that every positive
map from $M_m$ into $M_n$ is decomposable if and only if $(m,n)$ is
$(2,2)$, $(2,3)$ or $(3,2)$ \cite{{choi75},{stormer},{woronowicz}},
and many examples of indecomposable positive linear maps are known.
See the survey article \cite{kye_ritsu}.

Recall that a positive map is said to be exposed (respectively extremal) if it generates
an exposed (respectively extreme) ray in the convex cone of all positive linear maps.
Among various examples of indecomposable positive linear maps, exposed maps attract special interests in quantum
information theory. First of all, indecomposable exposed maps detect PPT entanglement of nonzero volume.
Furthermore, exposed maps are enough to determine separability, by the classical Straszewicz's Theorem:
Exposed points make a dense subset in the set of all extreme points. See \cite{rock}, Theorem 18.6.
In case of decomposable maps, the maps in (\ref{dec-exp}) are known to be exposed \cite{marcin_exp,yopp},
and there are no other exposed decomposable maps.

Nevertheless, there are very few known examples of exposed indecomposable positive linear maps in the literature.
Some variants \cite{cho-kye-lee} of the Choi map
\cite{choi-lam} between $M_3$ had been proved to be exposed in \cite{ha-kye_exp_Choi}. Note that the Choi map itself
is not exposed, even though it is extremal. On the other hand, the
Robertson's map \cite{robertson-83,robertson} between $M_4$ was shown \cite{chrus_Robert} to be
exposed. Some variants in the higher dimensional cases have been also considered
in \cite{sarbicki_exp}. As for linear maps between $M_2$ and $M_4$, a large class of
indecomposable exposed positive linear maps have been constructed in \cite{ha-kye_exp_woro},
motivated by the example of Woronowicz \cite{woronowicz-1}.

The author \cite{kye_3qb_EW} had considered multi-partite version of the duality,
and exhibited examples of indecomposable positive bi-linear maps with the {\sl full spanning property}
between $2\times 2$ matrices,
in order to get examples of three qubit separable states with full ranks.
A positive map with the full spanning property is automatically
indecomposable, and indecomposable exposed positive maps must have
the spanning property.
See the next section for the details.
Nevertheless, the spanning property does not imply the exposedness in general, even in the bi-partite systems. See the
example in Section 5 of \cite{ha-kye_uniq}.

The purpose of this note is to show that the positive bi-linear maps constructed in \cite{kye_3qb_EW}
are, in fact, exposed. The main tool is to use the four dimensional norm introduced in \cite{chen-han-kye}.
This norm turns out to be crucial to characterize \cite{chen-han-kye} separability of an important class of three qubit states,
like Greenberger-Horne-Zeilinger diagonal states. See also \cite{han_kye_GHZ,han_kye_phase}.
These provide the first concrete examples of indecomposable exposed positive {\sl bi-linear} maps.

In the next section, we review very briefly the duality between multi-partite separability
and positivity of multi-linear maps, and explain the motivation of the construction in \cite{kye_3qb_EW}.
In section 3, we compute the dual and the double dual to show that the map is exposed.

The author is grateful to Kyung-Hoon Han for useful discussion and comments.

\section{Duality and indecomposable positive bi-linear maps}

For a given multi-linear map $\phi:M_{d_1}\times \cdots \times M_{d_{n-1}}\to M_{d_n}$, we associate a matrix
$C_\phi\in M_d$ with $d=d_1 d_2 \cdots d_n$ by
$$
C_\phi
=\sum_{i_1,j_1,\dots, i_{n-1},j_{n-1}}
 |i_1\ran\lan j_1|\otimes \cdots\otimes |i_{n-1}\ran\lan j_{n-1}|\ot
 \phi(|i_1\ran\lan j_1|,\cdots, |i_{n-1}\ran\lan j_{n-1}|).
$$
The matrix $C_\phi$ is call the {\sl Choi matrix} of $\phi$.
The correspondence $\phi\mapsto C_\phi$ is just the Choi-Jamio\l kowski isomorphism \cite{choi75-10,jami}
for bi-partite case of $n=2$. See also \cite{CS-2009} for multi-partite cases.
For a state $\varrho\in M_d$ and an $(n-1)$-linear map $\phi$, we consider the bilinear pairing $\lan \varrho, \phi\ran$
defined by
$$
\lan\varrho,\phi\ran = \lan \varrho,C_\phi\ran = \tr (C_\phi \varrho^\ttt).
$$
Recall that a state in $M_d$ is said to be separable if it is the convex combination of pure product states, and
a multi-linear map $\phi:M_{d_1}\times \cdots \times M_{d_{n-1}}\to M_{d_n}$ is called positive if
$\phi(x_1,\dots x_{n-1})\in M_{d_n}$ is positive whenever $x_i\in M_{d_i}$ is positive for every $i=1,2,\dots,n-1$.
If $\varrho=|\xi\ran\lan\xi|$ is a pure product state with the product vector
$|\xi\ran=|\xi_1\ran\cdots |\xi_n\ran$ then we have the relation
$$
\lan\varrho,\phi\ran
=\lan C_\phi, |\xi\ran\lan\xi|\ran
=\lan \bar\xi |C_\phi|\bar\xi\ran
=\lan\bar\xi_n|\phi(|\xi_1\ran\lan\xi_1|,\dots, |\xi_{n-1}\ran\lan\xi_{n-1}|)|\bar\xi_n\ran.
$$
Therefore, we see \cite{kye_3qb_EW} that a multi-linear map $\phi$ is positive if and only if
$\lan\varrho,\phi\ran\ge 0$ for every separable state $\varrho$.
By the Hahn-Banach type separation theorem, we see that
$\varrho$ is separable if and only if $\lan\varrho,\phi\ran\ge 0$
for every positive multi-linear map $\phi$.
We also see that a self-adjoint matrix $W$ is the Choi matrix of a positive multi-linear map if and only if
$\lan \xi |W|\xi\ran\ge 0$ for every product vector $|\xi\ran$. Those self-adjoint matrices are called block-positive.
We say that a multi-linear map $\phi$ is {\sl completely positive} if its Choi matrix $C_\phi$ is positive
(semi-definite). See \cite{choi75-10} for linear cases.

For a given positive multi-linear map $\phi$, the dual $\phi^\prime$ of $\phi$ is defined by the set of all
(unnormalized) separable states $\varrho$ such that $\lan\varrho,\phi\ran=0$.
The set $\phi^\prime$ is an exposed face of the convex cone of all separable (unnormalized) states.
We denote by $P[\phi]$ the set of all product vectors $|\xi\ran$ such that $\lan\xi|C_\phi|\xi\ran=0$.
The double dual $\phi^{\prime\prime}$ is also defined by the set of all positive maps $\psi$ satisfying
$\lan\varrho,\psi\ran=0$ for every $\varrho\in\phi^\prime$. The double dual $\phi^{\prime\prime}$ is the
smallest exposed face containing the map $\phi$, and so we see that $\phi$ is exposed if and only if
$\phi^{\prime\prime}$ consists of nonnegative multiples of $\phi$.

For a given subset $S$ of $\{1,2,\dots,n\}$, the partial transpose  $A^{T(S)}$ is defined by
$$
(A_1\ot A_2\ot\cdots\ot A_n)^{T(S)}:=B_1\ot B_2\ot\cdots\ot B_n,
\quad \text{\rm with}\ B_j=\begin{cases} A_j^\ttt, &j\in S,\\ A_j,
&j\notin S,\end{cases}
$$
for multi-partite matrices.
A state $\varrho$ is said to be of PPT (positive partial transpose) if $\varrho^{T(S)}$ is positive for every subset
$S$. On the other hand, $\phi$ is called {\sl decomposable} if its Choi matrix
is the sum of partial transposes of positive matrices. It is clear that $\phi$ is decomposable if and only if
$\lan\varrho,\phi\ran\ge 0$ for every PPT state $\varrho$, and $\varrho$ is of PPT if and only if
$\lan\varrho,\phi\ran\ge 0$ for every decomposable map $\phi$.
On the other hand, the partial conjugate $\ket \xi^{\Gamma(S)}$ is defined by
$$
(\ket{\xi_1}\ot\cdots\ot\ket{\xi_n})^{\Gamma(S)}
:=\ket{\eta_1}\ot\cdots\ot\ket{\eta_n}, \quad \text{\rm with}\
\ket{\eta_j}=\begin{cases} \ket{\bar{\xi}_j}, &j\in S,\\
\ket{\xi_j}, &j\notin S,\end{cases}
$$
for product vectors.
We say that $\phi$ has the {\sl full spanning property} if
the set $\{|\xi\ran ^{\Gamma(S)}: \xi\in P[\phi]\}$ spans the whole space for every $S$.
It was shown in \cite{kye_3qb_EW} that the following are equivalent for a positive multi-linear map $\phi$:
\begin{itemize}
\item
$\phi$ has the full spanning property.
\item
The interior of the face $\phi^\prime$ lies in the interior of convex cone of all PPT states.
\item
The smallest exposed face $\phi^{\prime\prime}$ containing $\phi$ has no decomposable map.
\item
The set of PPT states $\varrho$ satisfying $\lan\varrho,\phi\ran < 0$
has a nonempty interior.
\end{itemize}
The last statement tells us that the map $\phi$ detects PPT entanglement with nonzero volume.
By the third condition, we also see that every indecomposable exposed map satisfies the above properties.

Now, we proceed to construct an example of indecomposable positive bi-linear map
$\phi:M_2\times M_2\to M_2$. It is a positive bi-linear map if and only if
the corresponding linear map
$M_2\to {\mathcal L}(M_2,M_2)$ sends positive matrices into positive maps
if and only if it sends positive matrices into decomposable maps by \cite{stormer}
if and only if the corresponding linear map
$M_2\to M_4$ sends positive matrices into matrices
which are sums of positive matrices and co-positive matrices, that is, partial transposes of positive matrices.
In other words, it is the sum of two maps $\phi_1$ and $\phi_2$, where
$\phi_1$ and $\phi_2$ send positive matrices into
positive and co-positive matrices, respectively.
If both of $\phi_1$ and $\phi_2$ are {\sl linear}  then
the corresponding bi-linear map $\phi$ must be decomposable.
Therefore, we have to look for {\sl non-linear} maps $\phi_1$ and $\phi_2$
with the above properties so that $\phi_1+\phi_2$ is {\sl linear}.
For every complex numbers, we consider the rank one matrix
$P_\alpha=\left(\begin{matrix} 1&\bar\alpha\\ \alpha&|\alpha|^2\end{matrix}\right)$, and put
$$
A_\alpha
=\left(\begin{matrix}
\cdot& \cdot &\cdot &\cdot\\
\cdot& |\bar\alpha+\alpha|^2 & \bar\alpha+\alpha &\cdot\\
\cdot &\alpha+\bar\alpha & 1& \cdot\\
\cdot& \cdot &\cdot &\cdot
\end{matrix}\right),\ \ \
B_\alpha
=\left(\begin{matrix}
\cdot& \cdot &\cdot &\cdot\\
\cdot& |\bar\alpha-\alpha|^2 & \bar\alpha-\alpha &\cdot\\
\cdot &\alpha-\bar\alpha & 1& \cdot\\
\cdot& \cdot &\cdot &\cdot
\end{matrix}\right),
$$
where $\cdot$ denotes zero.
Neither $P_\alpha\mapsto A_\alpha$ nor $P_\alpha\mapsto B_\alpha^T$ are linear.
But, we note that the map
$$
P_\alpha\mapsto A_\alpha + B_\alpha^T
=\left(\begin{matrix}
\cdot& \cdot &\cdot &\bar\alpha-\alpha\\
\cdot& 4|\alpha|^2 & \bar\alpha+\alpha &\cdot\\
\cdot &\bar\alpha+\alpha & 2& \cdot\\
\alpha-\bar\alpha& \cdot &\cdot &\cdot
\end{matrix}\right)
$$
is linear.

Motivated by the above discussion, we have considered in \cite{kye_3qb_EW}
the bilinear map $\phi:M_2\times M_2\to M_2$
which sends $([x_{ij}], [y_{ij}])\in M_2\times M_2$ to
\begin{equation}\label{def_phi}
\left(\begin{matrix}
sx_{22}y_{11}&
x_{12}y_{12}-x_{12}y_{21}+x_{21}y_{12}+x_{21}y_{21}\\
x_{12}y_{12}+x_{12}y_{21}-x_{21}y_{12}+x_{21}y_{21}
&t x_{11}y_{22}
\end{matrix}\right)\in M_2,
\end{equation}
where $s,t$ are positive numbers with $\sqrt{st}=2\sqrt 2$.
Note that the corresponding Choi matrix is given by
$$
C_\phi=\left(\begin{matrix}
\cdot&\cdot&\cdot&\cdot&\cdot&\cdot&\cdot &1\\
\cdot&\cdot&\cdot&\cdot&\cdot&\cdot&1&\cdot\\
\cdot&\cdot&\cdot&\cdot&\cdot&-1&\cdot&\cdot\\
\cdot&\cdot&\cdot&t&1&\cdot&\cdot&\cdot\\
\cdot&\cdot&\cdot&1&s&\cdot&\cdot&\cdot\\
\cdot&\cdot&-1&\cdot&\cdot&\cdot&\cdot&\cdot\\
\cdot&1&\cdot&\cdot&\cdot&\cdot&\cdot&\cdot\\
1&\cdot&\cdot&\cdot&\cdot&\cdot&\cdot&\cdot
\end{matrix}\right),
$$
if we use the lexicographic order
$000$, $001$, $010$, $011$, $100$, $101$, $110$, $111$.
It was shown \cite{kye_3qb_EW} that the map $\phi$ is indecomposable positive bi-linear map.
In fact, it was shown that $\phi$ has the full spanning property. See also \cite{han_kye_genuine}
for indecomposability.

We note that all the entries of $C_\phi$ are zero except for diagonal and anti-diagonals. Such matrices
are called {\sf X}-shaped, and are of the form
$$
X(a,b,c)= \left(
\begin{matrix}
a_1 &&&&&&& c_1\\
& a_2 &&&&& c_2 & \\
&& a_3 &&& c_3 &&\\
&&& a_4&c_4 &&&\\
&&& \bar c_4& b_4&&&\\
&& \bar c_3 &&& b_3 &&\\
& \bar c_2 &&&&& b_2 &\\
\bar c_1 &&&&&&& b_1
\end{matrix}
\right),
$$
for vectors
$a=(a_1,a_2,a_3,a_4),b=(b_1,b_2,b_3,b_4)\in\mathbb R^4$ and
$c=(c_1,c_2,c_3,c_4)\in\mathbb C^4$. Many important three qubit states arise in this form. For example,
Greenberger-Horne-Zeilinger diagonal states are {\sf X}-shaped, and $\varrho=X(a,b,c)$ is GHZ diagonal
if and only if $a=b$ and $c\in\mathbb R^4$.

\section{Exposedness}

We first find the dual face of $\phi$. To do this, we look for
pure product states $\varrho=|\xi\rangle\langle\xi|$ with the product vector
$|\xi\rangle=|x\rangle\otimes |y\rangle \otimes |z\rangle$ satisfying the relation
\begin{equation}\label{dual-vector}
0=\langle\varrho,\phi\ran=\lan C_\phi,|\xi\ran\lan\xi|\ran
=\lan\bar\xi |C_\phi|\bar\xi\ran=
\langle\bar z|\phi(|x\rangle\langle x|, |y\rangle\langle y|)|\bar z\rangle.
\end{equation}
We first consider the case when one of $|x\rangle$, $|y\rangle$ or $|z\rangle$ has a zero entry.
In this case, it is easy to see \cite{kye_3qb_EW}
that product vectors satisfying the relation  (\ref{dual-vector}) are one of the following:
\begin{equation}\label{pv_1}
\begin{aligned}
|x\ran 0\ran 1\ran=x_0|001\ran +x_1|101\ran,\quad
&|x\ran 1\ran 0\ran=x_0|010\ran +x_1|110\ran,\\
|0\ran y\ran 0\ran=y_0|000\ran +y_1|010\ran,\quad
&|1\ran y\ran 1\ran=y_0|101\ran +y_1|111\ran,\\
|0\ran 0\ran z\ran=z_0|000\ran +z_1|001\ran,\quad
&|1\ran 1\ran z\ran=z_0|110\ran +z_1|111\ran,
\end{aligned}
\end{equation}
with $|x\rangle, |y\rangle, |z\rangle\in\mathbb C^2$.

Next, we consider the case when $|\xi\ran$ has no zero entry.
In this case, the {\sf X}-part $\varrho_X$
of $\varrho=|\xi\ran\lan\xi|$ is again a separable state of rank four. In fact, it was shown in \cite{han_kye_phase}
that $\varrho_X$ is the average
$\varrho_X=\frac 14 \sum_{k=0}^3 |\xi(k)\ran\lan\xi(k)|$ of the four pure product states
given by
\begin{equation}\label{decom-rank1}
\begin{aligned}
|\xi(0)\ran=&|x_+ \ran \ot |y_+ \ran \ot |z_+ \ran,\\
|\xi(1)\ran=&|x_+ \ran \ot |y_- \ran \ot |z_- \ran,\\
|\xi(2)\ran=&|x_- \ran \ot |y_+ \ran \ot |z_- \ran,\\
|\xi(3)\ran=&|x_- \ran \ot |y_- \ran \ot |z_+ \ran,
\end{aligned}
\end{equation}
with $|x_\pm\ran =(x_0,\pm x_1)^\ttt$ and $|y_\pm\ran$, $|z_\pm\ran$
similarly. Furthermore, every rank four separable non-diagonal {\sf X}-state
arises in this way with a unique decomposition.

Because $C_\phi$ is {\sf X}-shaped, the relation (\ref{dual-vector}) is equivalent to
$\langle C_\phi,\varrho_X\rangle=0$.
Therefore, we proceed to find non-diagonal {\sf X}-state $\varrho=X(a,b,c)$ of rank four
so that $\langle\phi,\varrho\rangle=0$. Note \cite{han_kye_phase} that
$\varrho=X(a,b,c)$ is non-diagonal separable state of rank four if and only if
$$
a_ib_i=|c_j|^2\ (i,j=1,2,3,4),\qquad a_1a_4=a_2a_3,\qquad c_1c_4=c_2c_3.
$$
We may assume that $a_ib_i=|c_j|=1$.
Since $\lan C_\phi,\varrho\ran=ta_4+sb_4+2\re (c_1+c_2-c_3+c_4)$, we have
\begin{equation}\label{ineq1}
2\sqrt 2=\sqrt{st}\le \frac 12(ta_4+sb_4)=-\re (c_1+c_2-c_3+c_4).
\end{equation}
Putting $\bar c_1c_2=\bar c_3c_4=-e^{{\rm i}\theta}$, we also have
\begin{equation}\label{ineq2}
\begin{aligned}
2\sqrt 2\le
-\re (c_1+c_2-c_3+c_4)
&\le |c_1+c_2|+|c_3-c_4|\\
&=|1-e^{{\rm i}\theta}|+|1+e^{{\rm i}\theta}|
\le 2\sqrt 2.
\end{aligned}
\end{equation}
Therefore, all the inequality in (\ref{ineq1}) and (\ref{ineq2}) becomes the equalities to get
\begin{itemize}
\item
$ta_4=sb_4$ which implies $a_4=\sqrt{s/t}$,
\item
both $c_1+c_2$ and $-c_3+c_4$ are non-positive real numbers,
\item
$\theta=\frac \pi 2$ or $\theta=-\frac \pi 2$.
\end{itemize}
In the case of $\theta=\frac \pi 2$, we see that both $c_1+c_2=c_1(1-e^{{\rm i}\frac \pi 2})$
and $-c_3+c_4=-c_3(1+e^{{\rm i}\frac \pi 2})$ are non-positive real numbers, and so we have
$c=(\omega^{-3},\omega^3,\omega^{-1},\omega^{-3})$, where $\omega=e^{{\rm i}\frac \pi 4}$ is the
$8$-th root of unity. If $\theta=-\frac\pi 2$ then we also have
$c=(\omega^3,\omega^{-3},\omega,\omega^3)$. So, we have two kinds of rank four separable states
$$
\begin{aligned}
\varrho_1(a_1,a_2)
&=
X((a_1,a_2,ua_1a_2^{-1},u),(a_1^{-1},a_2^{-1},u^{-1}a_1^{-1}a_2,u^{-1}),(\omega^{-3},\omega^3,\omega^{-1},\omega^{-3})),\\
\varrho_2(a_1,a_2)
&=
X((a_1,a_2,ua_1a_2^{-1},u),(a_1^{-1},a_2^{-1},u^{-1}a_1^{-1}a_2,u^{-1}),(\omega^{3},\omega^{-3},\omega,\omega^{3}))
\end{aligned}
$$
in the dual $\phi^\prime$ of $\phi$,
with $u:=\sqrt{\frac st}$ and $a_1,a_2>0$.
This is enough to compute the double dual $\phi^{\prime\prime}$.

For the completeness, we find product vectors
$$
\begin{aligned}
|\xi\ran
&=(p_1,\alpha_1)^\ttt\ot (p_2,\alpha_2)^\ttt\ot (p_3,\alpha_3)^\ttt\\
&=(p_1p_2p_3, p_1p_2\alpha_3, p_1\alpha_2p_3, p_1\alpha_2\alpha_3,
\alpha_1p_2p_3, \alpha_1p_2\alpha_3, \alpha_1\alpha_2p_3, \alpha_1\alpha_2\alpha_3)^\ttt
\end{aligned}
$$
with $p_i>0$ and $\alpha_i\in\mathbb T$ satisfying the relation
$r\varrho_1=|\xi\ran\lan\xi|_X$ with $r>0$. Taking the {\sf X}-part, we have
$$
\begin{aligned}
r(a_1,a_2,ua_1a_2^{-1},u)
&=(p_1^2p_2^2p_3^2, p_1^2p_2^2, p_1^2p_3^2, p_1^2),\\
r(a_1^{-1},a_2^{-1},u^{-1}a_1^{-1}a_2,u^{-1})
&=(1,p_3^2,p_2^2,p_2^2p_3^2),\\
r(\omega^{-3},\omega^3,\omega^{-1},\omega^{-3})
&=p_1p_2p_3(\bar\alpha_1\bar\alpha_2\bar\alpha_3,\bar\alpha_1\bar\alpha_2\alpha_3,
\bar\alpha_1\alpha_2\bar\alpha_3,\bar\alpha_1\alpha_2\alpha_3),
\end{aligned}
$$
which implies $r=p_1p_2p_3=a_1$ together with
$$
\alpha_1^2=\omega^6,\quad \alpha_2^2=\omega^2,\quad \alpha_3^2=\omega^6,\quad
p_1^2=ua_1,\quad p_2^2=u^{-1}a_2,\quad p_3^2=a_1a_2^{-1}.
$$
Therefore, $(\alpha_1,\alpha_2)$ is one of $(\omega^3,\omega)$, $(\omega^3,\omega^5)$, $(\omega^7,\omega)$ and
$(\omega^7,\omega^5)$. In each case, $\alpha_3$ is automatically determined.
In this way, we have the following four kinds
\begin{equation}\label{pv_2}
\begin{aligned}
\eta_1(a_1,a_2)
&=(u^{\frac 12}a_1^{\frac 12},\omega^3)^\ttt\ot (u^{-\frac 12}a_2^{\frac 12},\omega)^\ttt\ot (a_1^{\frac 12}a_2^{-\frac 12},\omega^7)^\ttt,\\
\eta_2(a_1,a_2)
&=(u^{\frac 12}a_1^{\frac 12},\omega^3)^\ttt\ot (u^{-\frac 12}a_2^{\frac 12},\omega^5)^\ttt\ot (a_1^{\frac 12}a_2^{-\frac 12},\omega^3)^\ttt,\\
\eta_3(a_1,a_2)
&=(u^{\frac 12}a_1^{\frac 12},\omega^7)^\ttt\ot (u^{-\frac 12}a_2^{\frac 12},\omega)^\ttt\ot (a_1^{\frac 12}a_2^{-\frac 12},\omega^3)^\ttt,\\
\eta_4(a_1,a_2)
&=(u^{\frac 12}a_1^{\frac 12},\omega^7)^\ttt\ot (u^{-\frac 12}a_2^{\frac 12},\omega^5)^\ttt\ot (a_1^{\frac 12}a_2^{-\frac 12},\omega^7)^\ttt
\end{aligned}
\end{equation}
of product vectors satisfying (\ref{dual-vector}),
up to constant multiples. In case of $\theta=-\frac\pi 2$, we consider the relation
$r\varrho_2=|\xi\ran\lan\xi|_X$ to get the following four product vectors
\begin{equation}\label{pv_3}
\begin{aligned}
\zeta_1(a_1,a_2)
&=(u^{\frac 12}a_1^{\frac 12},\omega^5)^\ttt\ot (u^{-\frac 12}a_2^{\frac 12},\omega^7)^\ttt\ot (a_1^{\frac 12}a_2^{-\frac 12},\omega)^\ttt,\\
\zeta_2(a_1,a_2)
&=(u^{\frac 12}a_1^{\frac 12},\omega^5)^\ttt\ot (u^{-\frac 12}a_2^{\frac 12},\omega^3)^\ttt\ot (a_1^{\frac 12}a_2^{-\frac 12},\omega^5)^\ttt,\\
\zeta_3(a_1,a_2)
&=(u^{\frac 12}a_1^{\frac 12},\omega)^\ttt\ot (u^{-\frac 12}a_2^{\frac 12},\omega^7)^\ttt\ot (a_1^{\frac 12}a_2^{-\frac 12},\omega^5)^\ttt,\\
\zeta_4(a_1,a_2)
&=(u^{\frac 12}a_1^{\frac 12},\omega)^\ttt\ot (u^{-\frac 12}a_2^{\frac 12},\omega^3)^\ttt\ot (a_1^{\frac 12}a_2^{-\frac 12},\omega)^\ttt.\\
\end{aligned}
\end{equation}
In total, we found all the product vectors satisfying (\ref{dual-vector})
listed in (\ref{pv_1}), (\ref{pv_2}) and (\ref{pv_3}).
We note that four product vectors $\eta_1,\eta_4,\zeta_1,\zeta_4$ were missing in \cite{kye_3qb_EW}, and correct here.

Now, we proceed to find the double dual dual $\phi^{\prime\prime}$ of $\phi$. In other word, we are looking for
self-adjoint block-positive three qubit matrices $W$ satisfying the relation $\lan W,\varrho\ran=0$
for every $\varrho\in\phi^\prime$.
We see from the list (\ref{pv_1}) that $\phi^\prime$ contains the pure
product state associated with the following product vectors:
\begin{equation}\label{pv-4}
|0\ran 0\ran 0\ran,\quad
|0\ran 0\ran 1\ran,\quad
|0\ran 1\ran 0\ran,\quad
|1\ran 0\ran 1\ran,\quad
|1\ran 1\ran 0\ran,\quad
|1\ran 1\ran 1\ran.
\end{equation}
Therefore, the corresponding diagonal entries of $W$ must be zero. From this, we show that
$W$ must be {\sf X}-shaped.
We note that the $4\times 4$ principle submatrix $W_0$ of $W$ with entries from $\{000, 001, 010, 011\}$ is block-positive
with zero diagonals for $000,001,010$ entries. Because the $2\times 2$ principle submatrices with entries
$\{000,001\}$, $\{000, 010\}$ $\{010,011\}$, $\{001,011\}$ are positive, we see that $W_0$ is of the form
$$
W_0=
\left(\begin{matrix}
\cdot &\cdot &\cdot &u\\
\cdot &\cdot &v &\cdot\\
\cdot &\bar v &\cdot &\cdot\\
\bar u &\cdot &\cdot & 1
\end{matrix}\right)
$$
up to scalar multiplication.
We note that the corresponding positive map
sends the positive matrix $P_\alpha:=\left(\begin{matrix}1&\alpha\\ \bar\alpha
&|\alpha|^2\end{matrix}\right)$ to the $2\times 2$  matrix
$\left(\begin{matrix}
0& \alpha u+\bar\alpha\bar v\\
\bar\alpha \bar u+ \alpha v & |\alpha|^2
\end{matrix}\right)$
which must be positive. This implies that $\alpha u+\bar\alpha\bar v=0$ for every complex $\alpha$,
from which we conclude that $u=v=0$. Therefore, all the entries of $W_0$ are zero except for $011$ place.
In the above argument, we fix $0$ in the first place of indices to get $4\times 4$ principle submatrix
which is still block-positive. We fix $0$ or $1$ in the first, second or third places to get six $4\times 4$
principle submatrices which are still block-positive. All the entries of such principle submatrices are zero
except for diagonal entries at $011$ and $100$ places. This implies that $W$ must be {\sf X}-shaped.

Now, we can write $W=X((0,0,0,x_4),(0,0,0,y_4),(z_1,z_2,z_3,z_4))$. By \cite{han_kye_tri}, we know that
$W$ is block-positive if and only if
$$
\sqrt {x_4y_4}\ge
\max_\theta\left(|z_1 e^{{\rm i}\theta}+\bar z_4| + |z_2
e^{{\rm i}\theta} +\bar z_3|\right):=\|z\|_\xx.
$$
See \cite{chen-han-kye} for properties of the norm $\|\ \|_\xx$ defined as above.
Considering the state
$$
\varrho
:=\frac 12(\varrho_1(1,1)+\varrho_2(1,1))\\
=X((1,1,u,u),(1,1,u^{-1},u^{-1}),{\textstyle \frac 1{\sqrt 2}}(-1,-1,1,-1))
$$
in the dual $\phi^\prime$, we have
$$
0=\lan \varrho,W\ran=x_4u+y_4u^{-1}-\sqrt 2\re (z_1+z_2-z_3+z_4).
$$
Therefore, we have
\begin{equation}\label{key-ineq}
\begin{aligned}
2\sqrt 2\|z\|_\xx
&\le\sqrt 2\cdot 2\sqrt{x_4y_4}\\
&\le \sqrt 2(x_4\sqrt{s/t}+y_4\sqrt{t/s})\\
&=2\re (z_1+z_2-z_3+z_4)\\
&\le 2(|z_1|+|z_2|+|z_3|+|z_4|)\le 2\sqrt 2\|z\|_\xx.
\end{aligned}
\end{equation}
See \cite{chen-han-kye} and Appendix of \cite{han_kye_GHZ}
for the last inequality $\|z\|_\xx\ge \frac 1{\sqrt 2} \|z\|_1$, where the equality holds
only if the entries have a common magnitude and
$$
(\arg z_1+\arg z_4)-(\arg z_2+\arg z_3)=\pi
$$
by Proposition 3.8 and (27) in \cite{chen-han-kye}.
Therefore, we have the following:
\begin{enumerate}
\item[(i)]
$\sqrt 2\|z\|_\xx= \|z\|_1$, and so $z$ has a common magnitude,
\item[(ii)]
$x_4\sqrt{s/t}=y_4\sqrt{t/s} \ \Longleftrightarrow\ x_4s=y_4t$,
\item[(iii)]
$\re (z_1+z_2-z_3+z_4)=|z_1|+|z_2|+|z_3|+|z_4|$.
\end{enumerate}
From (i) and (iii), we see that $z=r_1(1,1,-1,1)$ for a nonnegative $r_1>0$, and we also have $(x_4,y_4)=r_2(t,s)$
for $r_2>0$ by (ii). Then we have $2\re (z_1+z_2-z_3+z_4)=8r_1$ and
$\sqrt 2(x_4\sqrt{s/t}+y_4\sqrt{t/s})=\sqrt 2\cdot 2\sqrt{st} \cdot r_2=8r_2$. Therefore, we have $r_1=r_2$
by the equality in (\ref{key-ineq}), and finally
conclude that $W$ is a nonnegative multiple of $C_\phi$. This completes the proof that the bi-linear map
$\phi$ is exposed.


\end{document}